\documentclass[11pt,amsmath,amssym,latexsymb]{article}
\usepackage{amssymb,latexsym}
\usepackage{longtable,tabularx}
\usepackage{makeidx}
\usepackage{graphicx}
\usepackage{psfrag}
\begin {document}

\def\hat{\widehat}
\def\tilde{\widetilde}
\def \bfo {\begin {eqnarray*} }
\def \efo {\end {eqnarray*} }
\def \ba {\begin {eqnarray*} }
\def \ea {\end {eqnarray*} }
\def \beq {\begin {eqnarray}}
\def \eeq {\end {eqnarray}}
\newtheorem {proposition}{Proposition}
\newtheorem {theorem}{Theorem}

\def\picture #1 by #2 (#3){
      \vsquare to #2{
        \hrule width #1 height 0pt depth 0pt
        \hfill
        \special{picture #3}}}

\def\scaledpicture #1 by #2 (#3 scaled #4){{
      \dimen0=#1 \dimen1=#2
      \divide\dimen0 by 1000 \multiply\dimen0 by #4
      \divide\dimen1 by 1000 \multiply\dimen1 by #4
      \picture \dimen0 by \dimen1 (#3 scaled #4)}}

\def \R {{\bf {R}}}

\def \H2s {H^{s+1}_0(\partial M\times [0,T/2])}

\def \det {\hbox{det}}

\def \l {\sigma}

\def \pa0 {\partial _0}

\def \p {\partial}

\def\M{{M}}

\def \hat {\widehat }


\def\tilde{\widetilde}


\def \mbeq {\begin {eqnarray}}
\def \meeq {\end {eqnarray}}


\def \bfo {\begin {displaymath} }
\def \efo {\end {displaymath} }

\def \beq {\begin {eqnarray}}
\def \eeq {\end {eqnarray}}
\def \ba {\begin {eqnarray*}}
\def \ea {\end  {eqnarray*}}

\def \R {{\bf {R}}}

\def \H2s {H^{s+1}_0(\partial \M\times [0,T/2])}

\def \det {\hbox{det}}

\def \l {\sigma}

\def \pa0 {\partial _0}

\def \p {\partial}



\def \tilde{\widetilde}

\title{On nonuniqueness for Calder\'on's inverse problem}

\author{
Allan Greenleaf, Matti Lassas and Gunther Uhlmann
\footnote{University of Rochester, Rochester, NY 14618; Rolf Nevanlinna
Institute, University of Helsinki Helsinki, P.O.Box 4, FIN-00014, Finland;
and University of Washington, Seattle, WA 98195
} }
\date{}
\maketitle
{\bf Abstract.}
We construct anisotropic conductivities with  the same
Dirichlet-to-Neumann map as a homogeneous isotropic conductivity. These
conductivities are  singular close to a surface
inside the body.
\bigskip

\section{Introduction}

    An
anisotropic conductivity on a domain $\Omega\subset\R^n$ is defined by a
symmetric, positive
semi-definite matrix-valued function,
$\l=(\l^{ij}(x))$. In the absence of sources or sinks, an electrical
potential $u$ satisfies
\beq
(\nabla\cdotp \sigma\nabla)u= \p
_j\sigma^{jk}(x) \p_k u &=& 0\hbox{ in } \Omega,\\
u|_{\p \Omega}&=&f,\nonumber
\eeq
where $f$ is the prescribed voltage on the boundary.
Above and hereafter we use the Einstein summation convention
where there is no danger of confusion. The resulting
voltage-to-current (or
Dirichlet-to-Neumann) map is then defined by
\beq
\Lambda_\sigma(f)= Bu|_{\p \Omega}
\eeq
where
\beq
Bu= \nu_j\sigma^{jk}\p_k u,
\eeq
$u$ is the solution of (1) and $\nu=(\nu_1,\dots,\nu_n)$ is the unit
normal vector of $\p\Omega$.

Applying
the divergence theorem, we have
\beq
Q_\l (f)=:\int_\Omega  \sigma^{jk}(x)\frac {\p u}{\p x^j}
\frac {\p u}{\p x^k}
dx=\int_{\p \Omega} \Lambda_\l(f) f dS,
\eeq
where $u$ solves (1) and $dS$ denotes surface measure on $\p \Omega$.
$Q_\l (f)$ represents the
power needed to maintain the potential $f$ on $\p \Omega$. By (4),
knowing $Q_\l$ is
equivalent with knowing
$\Lambda_\l$. If $F:\Omega\to \Omega,\quad F=(F^1,\dots,F^n)$, is a
diffeomorphism with $F|_{\p
\Omega}=\hbox{Identity}$, then by making the change of variables $y=F(x)$
and setting $u=v\circ F^{-1}$ in the first
integral in (4), we obtain
\ba
\Lambda_{F_*\l}=\Lambda_\l,
\ea
where
\beq\label{cond and metr}
(F_*\sigma)^{jk}(y)=\left.
\frac 1{\det [\frac {\p F^j}{\p x^k}(x)]}
\sum_{p,q=1}^n \frac {\p F^j}{\p x^p}(x)
\,\frac {\p F^k}{\p x^q}(x)  \sigma^{pq}(x)\right|_{x=F^{-1}(y)}
\eeq
is the push-forward of the conductivity $\l$ by $F$. Thus, there is a large
(infinite-dimensional) class of conductivities which give rise to the
same electrical
measurements at the boundary. The uniqueness question, first proposed by
Calder\'on \cite{C},  we wish to
address is whether two
conductivities with the same Dirichlet-to-Neumann map must be such
pushforwards of each
other. By a direct construction, we will show that the answer is no.
Let  $D\subset\subset  \Omega$
be a smooth subdomain.
     We will construct in
$\Omega$ a conductivity $\hat \sigma$ for which
the boundary measurements coincide with
those made for homogeneous  conductivity $\gamma=1$ in $\Omega$.
We note that the  conductivity
    is singular in the sense that some components of
the conductivity tensor  go to zero, i.e., correspond to
perfectly insulating directions  and some
components can go to infinity, i.e., correspond to
perfectly conducting directions, as one approaches $\partial D$.
The physical meaning of these counterexamples is evident:
For instance in medical imaging, e.g., in Electrical Impedance
Tomography (EIT) certain anisotropic structures can form
barriers which, including their interior, appear in measurements
to be a homogeneous medium.

To construct these counterexamples, we need to consider a variant of
(1), which comes from
the Laplace-Beltrami operator on a compact Riemannian manifold with
boundary. Let us assume now that $(M,g)$ is an $n$-dimensional
Riemannian manifold with smooth
boundary $\p M$. The metric
$g$ is assumed to be symmetric and positive definite. The invariant
object analogous to the
conductivity equation (1) is the Laplace-Beltrami operator, which is given by
\beq
\Delta_g u=G^{-1/2}\p_j( G^{1/2}g^{jk}\p_k u)
\eeq
where $G=\det (g_{jk}),$ $[g_{jk}]=[g^{jk}]^{-1}$. The
Dirichlet-to-Neumann map is defined
by solving the Dirichlet problem
\beq\label{Laplace 3}
& &\Delta_g u=0\quad\hbox{in}\quad M,\\
& &u|_{\p  M}=f. \nonumber
\eeq
The operator analogous to $\Lambda_\l$ is then
\beq
\Lambda_g(f)=G^{1/2} \nu_j g^{jk}\frac{\p u}{\p x_k}|_{\p M},
\eeq
with $\nu=(\nu_1,\dots,\nu_n)$ the outward unit normal to $\p M$. In
dimension three or
higher, the conductivity matrix and the Riemannian metric are related by
\beq\label{sigma-g}
\sigma^{jk}=\det(g)^{1/2}g^{jk},\quad\hbox{or}\quad
g^{jk}=\det(\sigma)^{2/(n-2)}\sigma^{jk}.
\eeq
Moreover, $\Lambda_g=\Lambda_\l$, and $\Lambda_{\psi^*g}=\Lambda_g$, where
$\psi^*g$ denotes the pullback of the metric $g$ by a diffeomorphism
of $M$ fixing $\p M$
\cite{LeU}.

In dimension two, (9) is not valid; in this case, the conductivity
equation can be
reformulated as
\beq\label{Laplace}
& &\hbox{Div}_g (\beta\, \hbox{Grad}_g u)=0\quad\hbox{in}\quad {M},\\
& &u|_{\p  {M}}=f \nonumber
\eeq
where $\beta$ is the scalar function  $\beta=|\det\, \sigma|^{1/2}$,
$g=(g_{jk})$ is equal to $(\sigma_{jk})$,
and $\hbox{Div}_g$ and $\hbox{Grad}_g$ are the divergence and
gradient operators with
respect to the Riemannian metric $g$. Thus we see that, in two
dimensions, Laplace-Beltrami
operators correspond only to those conductivity equations for which
$\det(\sigma)=1$.

For domains in two dimensions, Sylvester\cite{Sy} showed, using
isothermal coordinates, that one can reduce the anisotropic problem
to the isotropic one;
combining this with the isotropic result of Nachman\cite{Na}, one obtains

\begin{theorem} If $\sigma$ and $\widetilde \sigma$
are two $C^3$-smooth anisotropic conductivities in $\Omega\subset \R^2$
    for which $\Lambda_\sigma=\Lambda_{\widetilde\sigma}$, then there
is a diffeomorphism $F:\Omega\to \Omega$,
$F|_{\p \Omega}=Id$ such that
$
\widetilde \sigma=F_*\sigma.
$
\end{theorem}

In dimensions three and higher, the following
result is known (see \cite{LU}, \cite{LTU}, and  \cite{LeU}):

\begin{theorem}
If $n\ge 3$ and $(M,\partial M)$ is a $C^\omega$ manifold with
connected, $C^\omega$
boundary, and $g,\tilde{g}$ are $C^\omega$ metrics on $M$ such that
$\Lambda_g=\Lambda_{\tilde{g}}$, then there exists a $C^\omega$
diffeomorphism $F:M\to M$
such that $F|_{\p D}=Id$.
\end{theorem}

\section{Counterexamples}

Returning now to domains $\Omega\subset\subset\R^n, n\ge 3$, let
$D\subset\subset \Omega$ be an open subset with smooth boundary
and $g=g_{ij}$ be a metric on
$\Omega$. Let $y\in D$ be such that there is
a diffeomorphism $F:\Omega\setminus\{y\}\to
\Omega\setminus \overline D$, and let $\tilde g=F_*g$ on
$\Omega\setminus \overline D$. To
obtain a conductivity on all of $\Omega$, first  extend
$\tilde g$ to a bounded metric
inside $D$ and denote this new metric on $\Omega$ by $\hat g$.
We make this continuation so that the conductivity jumps
on $\p D$ and that $\hat g\geq c>0$ is smooth inside $D$, but is
otherwise arbitrary.
Let $\hat \sigma$ be the  conductivity corresponding to $\hat g$ by
(9) in $\Omega$.
We say that $v$ is a solution of the conductivity equation if
\beq\label{main eq}
& &\nabla\cdotp \hat\sigma\nabla v(x)=0\hbox{ in the sense of distributions
in }\Omega,\\
& &v|_{\p \Omega}=f_0, \nonumber \\
& &v\in L^\infty(\Omega) \nonumber
\eeq
where $f_0\in H^\frac12(\p \Omega)$.
That the equation in the sense of distributions
means that $v\in H^1(\Omega)$ and $\sigma\nabla v\in
H(\Omega;\nabla\cdotp)
=\{w\in L^2(\Omega):\
\nabla\cdotp w\in L^2(\Omega)\}$.
If this problem has a
unique solution, we define the Dirichlet-to-Neumann map
\ba
\Lambda_{\hat \sigma}f_0= n\cdotp \hat\sigma \nabla v|_{\p \Omega}.
\ea
Our aim is to prove the following
result, first announced in \cite{GLU}:

\begin{theorem}\label{main}
Let $\Omega\subset \R^n$, $n\geq 3$, and $g=g_{ij}$ a metric on $\Omega$.
Let $D\subset \Omega$ be such there is
a $C^\infty$-diffeomorphism $F:\Omega\setminus\{y\}\to
\Omega\setminus \overline D$ satisfying
$F|_{\p \Omega}=Id$ and that
\beq\label{Q 4}
dF(x)\geq c_0I,\quad
\hbox{det}\,(dF(x))\geq c_1\,\hbox{dist}_{{}_{\R^n}}\,(x,y)^{-1}
\eeq
where $dF$ is the Jacobian matrix in Euclidean coordinates
of $\R^n$ and $c_0,c_1>0$.
Let $\tilde g=F_*g$ and $\hat g$ be an extension
of $\tilde g$ into $D$ such that it is positive definite
in $D^{int}$. Finally, let  $\gamma$ and $\hat \sigma$
be the conductivities corresponding to $g$ and $\hat g$
by (9).
Then the boundary value problem for the conductivity equation with conductivity
$\hat \sigma$ is uniquely solvable and
\ba
\Lambda_{\hat \sigma}=\Lambda_\gamma.
\ea
\end{theorem}

\noindent
Note that here is no diffeomorphism $H:\Omega\to \Omega$ such that
$\hat \sigma=H_*\gamma$, so the Riemannian manifolds corresponding
to $\hat \sigma$ and $\gamma$ cannot be the same.
Also, $\hat \sigma$ can be changed in arbitrary way inside $D$
without changing boundary measurements.
The proof of Theorem \ref{main}  will be given below.
\medskip

\noindent
{\bf Example.} Let $\Omega=B(0,2)\subset \R^3$ be the ball with center
0 and radius 2.
Consider $y=0\in D=B(0,1)$ and the map
$F:\Omega\setminus\{0\}\to
\Omega
\setminus \overline D$ given by
\beq\label{A 2}
F:x\mapsto (\frac {|x|}2+1)\frac x{|x|}.
\eeq
Let $\gamma=1$ be the homogeneous conductivity in $\Omega$
and define $\sigma=F_*\gamma$.
Now the metric tensor $g$ and
the corresponding conductivity $\sigma_g$ are related by
$\sigma_g=|\hbox{det}\,g|^{1/2}g^{jk}$.
Let $g$ be the metric corresponding to $\gamma $
and $\tilde g$ be the metric corresponding to $\sigma$.
Consider these  in
the standard spherical coordinates on $\Omega\setminus \{0\}$,
   $(r,\phi,\theta)\mapsto
     (r\sin\theta \cos \phi,r\sin\theta \sin \phi,r\cos\theta)\in \R^3$.
With respect to these coordinates, we see that the metric
$g$ and conductivity $\gamma$ correspond to the matrices
\ba& &
g=\left(\begin{array}{ccc}
1 & 0 & 0\\
0 & r^2& 0 \\
0 & 0 &r^2\sin^2 \theta\\
\end{array}
\right), \quad \gamma=\left(\begin{array}{ccc}
r^2\sin \theta & 0 & 0\\
0 & \sin \theta& 0 \\
0 & 0 &(\sin \theta)^{-1}\\
\end{array}
\right)
\ea
and  $\tilde g$ and $\sigma$ correspond in the annulus $1<r<2$ to the matrices
\ba
& & \tilde g=
\left(\begin{array}{ccc}
4 & 0 & 0\\
0 & 4(r-1)^2 & 0 \\
0 & 0 &   4(r-1)^2 \sin^2 \theta\\
\end{array}
\right),\\
& &   \sigma=
\left(\begin{array}{ccc}
2(r-1)^2\sin \theta & 0 & 0\\
0 & 2 \sin \theta & 0 \\
0 & 0 &  2 (\sin \theta)^{-1}\\
\end{array}
\right).
\ea
Let $\hat\sigma$ be a continuation of $\sigma$ that is
     $C^\infty$-smooth
in $D$. Then these metric are as in
for Theorem \ref{main} and in particular $F$ satisfies
(\ref{Q 4}).
\medskip

To start to prove Theorem \ref{main}, we start with the following result:

\begin{proposition}
Let $\Omega\subset \R^n$, $n\geq 3$, and $g=g_{ij}$ a metric on $\Omega$. Let
$u$ satisfy
\ba
& &\Delta_{g} u(x)=0\quad\hbox{in }\Omega,\\
& &u|_{\p \Omega}=f_0\in C^\infty(\p \Omega).
\ea
Let $D\subset \Omega$ be  such that there is
a diffeomorphism $F:\Omega\setminus\{y\}\to
\Omega\setminus \overline D$ satisfying
$F|_{\p \Omega}=Id$. Let $\tilde g=F_*g$ and $v$ be a function
satisfying
\ba
& &\Delta_{\tilde g} v(x)=0\quad\hbox{in }\Omega\setminus \overline D,\\
& &u|_{\p \Omega}=f_0,\\
& &u\in L^\infty(\Omega\setminus \overline D).
\ea
Then $u$ and $F^*v$ coincide and  have the same Cauchy data on
$\p\Omega$,
\beq\label{C 2}
\p_\nu u|_{\p M}=\p_{\tilde \nu} F^*v|_{\p M}
\eeq
where $\nu$ is unit normal vector in metric $g$ and
$\tilde \nu$ is unit normal vector in metric $\tilde g$.
Moreover, for the constant $c_0:=u(y)$ we have
\ba
\lim_{x\to \p D}v(x)=c_0.
\ea
\end{proposition}
{\bf Proof.}  Let $g=(g_{ij})$ be a
Riemannian metric tensor defined on  $\Omega\subset \R^n$.

First, we continue the $C^\infty$-metric $g$
to a metric ${\frak g}$ in $\R^n$ such that ${\frak g}_{ij}(x)=
     g_{ij}(x)$ for $x\in \Omega$ is such that for  any $y\in \R^m$
there is a positive Green's function $G(x,y)$
satisfying
\ba
-\Delta_{\frak g} G(\cdotp,y)=\delta_y\quad\hbox{in }\R^n.
\ea
There are several easy ways to obtain this continuation.
For instance, we can continue $g$ to a metric ${\frak g}$
such that outside some ball
$B(0,R)$ the metric is hyperbolic.
This implies that the manifold $(\R^n,{\frak g})$
is non-parabolic and has a positive non-constant
super-harmonic function (see  \cite{Gr}).
By \cite{LiT} there thus exists
a positive Green's function $G(x,y)$.

Next we consider the probability that Brownian
motion $B_t^x$ on manifold $(\R^n,{\frak g})$
sent from the point $x$ at time $t=0$
enters an open  set $U$: Let
\ba
e_U(x)=P(\{\hbox{there is $t>0$ such that $B^x_t\in U$}\})
\ea
Let $A=\R^n\setminus U$. By Hunt's theorem (see,
\cite{H}, or \cite{Gr}, Prop. 4.4)
we have
\ba
e_U(x)=s_A(x)
\ea
where $s_A(x)$ is the super-harmonic potential of $A$, that is,
   $s_A(x)$ is infimum of all bounded super-harmonic functions
$h$ in $\R^n$ such that $h|_U=1$ and $h\geq 0$ in $\R^n$.

Let $U=B(y,r)$ and
\ba
m(y,r)=\inf_{x\in \p B(y,r)} G(x,y).
\ea
The function
\ba
h_{y,r}(x)=\frac 1{m(y,r)}\min \big(G(x,y),m(y,r)\big)
\ea
is positive super-harmonic function and satisfies $h_{y,r}|_{B(y,r)}=1$.
Since $e_U(x)\leq h_{y,r}(x)$ and
$
\lim_{r\to 0}m(y,r)=\infty $, we see that
\beq\label{Prob of hitting}
\lim_{r\to 0} P(\{\hbox{there is $t>0$ such that $B^x_t\in B(y,r)
$}\})=0.
\eeq
In particular, taking limit $r\to 0$ we see that the probability
that $B^x_t=y$ for some $t>0$ is zero.

Now,  let $\tilde u$ be any solution of
\beq\label{C 1}
& &\Delta_{g}\tilde u(x)=0\quad\hbox{in }\Omega\setminus \{y\}\\
& &\tilde u|_{\p \Omega}=f_0, \nonumber  \\
& &\tilde u\in L^\infty(\Omega\setminus \{y\}).\nonumber
\eeq
Denote $f_r=\tilde u|_{\p B(y,r)}$
and let
$\tau(r,x)\in (0,\infty]$ be the first time when
$B^x_t\in \p \Omega\cup \p B(y,r)$. Similarly,
let $\tau(0,x)$ be the hitting time to $\p \Omega$.

Let next $\chi_{0,r}^x$ be a random variable defined by
$\chi^x_{0,r}=1$ if $B_{\tau(r,x)}^x\in \p \Omega$
and zero otherwise, and let $\chi^x_{1,r}=
1-\chi^x_{0,r}$.
Then by Kakutani's formula,
\ba
\tilde u(x)={\Bbb E}(\chi^x_{0,r}f_0(B_{\tau(r,x)}^x))+{\Bbb
E}(\chi^x_{1,r}f_r(B_{\tau(r,x)}^x)).
\ea
Letting $r\to 0$ and using the fact that $||f_r||_\infty\leq
||\tilde u||_\infty$
are uniformly bounded and (\ref{Prob of hitting}) we see
$\tilde u(x)={\Bbb E}(f_0(B_{\tau(0,x)}^x))$. Thus
   $\tilde u(x)=u(x)$ for $x\in
\Omega\setminus \{y\}$ where
\beq\label{Q 2}
& &\Delta_g u(x)=0\quad\hbox{in }\Omega,\\
& &u|_{\p \Omega}=f_0.\nonumber
\eeq

Thus we have shown that
boundary value problem
(\ref{C 1}) is solvable, and that the solution is unique.
Next we change this problem
to an equivalent problem.
Let $F:\Omega\setminus\{y\}\to
\Omega\setminus \overline D$
be a diffeomorphism
     that is identity at boundary $\p \Omega$.
Define the metric $\tilde g=F_*g$; then the boundary value problem (\ref{Q 2})
is equivalent to
\beq\label{Q 3}
& &\Delta_{\tilde g} v(x)=0\quad\hbox{in }\Omega\setminus \overline D\\
& &v|_{\p \Omega}=f_0, \nonumber  \\
& &v\in L^\infty(\Omega\setminus \overline D).\nonumber
\eeq
and solutions of the problems  (\ref{Q 2}) and (\ref{Q 3})
are related by $v(x)=u(F(x))$ for $x\in \Omega\setminus \{y\}$.

Clearly, the Cauchy data of the equations (\ref{Q 2})
and  (\ref{Q 3}) coincide in the sense of (\ref{C 2}).
Moreover,
\ba
\lim_{x\to \p D}v(x)=c_0:=u(y)
\ea
This concludes the proof of Prop. 1.  $\square$

\medskip
We  remark that the application of Brownian motion above is not
essential but
makes the proof  perhaps more intuitive.
 Alternatively, in the proof of Prop. 1 one can use properties of
$L^p$-Sobolev spaces. Indeed, for $\tilde u$ satisfying
$
\Delta_{g}\tilde u(x)=0$ in $\Omega\setminus \{y\}$, $\tilde u|_{\p \Omega}=f_0$,
and $\tilde u\in L^\infty(\Omega\setminus \{y\})$,
we can consider the extension $\tilde u\in L^\infty(\Omega)$ that satisfies
\ba
\Delta_{g}\tilde u(x)=F\quad\hbox{in }\Omega,
\ea
where $F$ is a distribution
is supported in $y$.
Then $F$  has to be a finite sum of derivatives of the 
Dirac delta distribution 
supported at $y$. Now, $F\in W^{-2,p}(\Omega)$ for all $1<p<\infty$,
and since $\delta_y\not \in W^{-2,p}(\Omega)$ for $p>\frac n{n-2}$
and the same is true for the derivatives of the delta distribution, 
we see that that $F=0$. This implies
that $\tilde u=u$. 
\medskip

We now turn to the proof of Theorem 3. From now on, we assume that
$F:\Omega\setminus\{y\}\to
\Omega\setminus \overline D$, $F(x)=(F^1(x),\dots,F^n(x))$ is such
that condition (\ref {Q 4}) is satisfied. First,  continue  $\tilde
g$ into $D$ so that
that it is positive definite
in $D^{int}$. Next, extend $v$ inside $D$ to a function
\ba
& &h(x)=v(x)\hbox{ for }x\in \Omega\setminus \overline D,\\
& &h(x)=c_0\hbox{ for }x\in \overline D.
\ea
   Our aim is to show that $h$ is the solution of
(\ref{main eq}). Since any solution of (\ref{Q 3}) is constant
on the boundary of $\p D$ and $\tilde g$ is positive definite
in $D$, we see that the solution has to be constant inside $D$.
Thus $h$ is the unique solution if it is a solution.

Now we are ready to show that $v$ is a solution also in the sense of
distributions.
First we note that when $y(x)=F^{-1}(x)$ we have $v(x)=u(F^{-1}(x))$
and
\ba
\frac {\p v}{\p x^j}(x)=\frac {\p u}{\p y^k}(y(x))\frac {\p y^k}{\p
x^j}(x)
\ea
and since $u\in H^1(\Omega)$ and $\frac {\p y^k}{\p x^j}\in
L^\infty(\Omega\setminus \overline D)$ we have $v\in
H^1(\Omega\setminus \overline D)$.
Also, as $h\in C(\Omega\setminus D^{int})$
and the trace $v\mapsto v|_{\partial D}$
is a continuous map $C(\Omega\setminus D^{int})\cap
    H^1(\Omega\setminus \overline D)\to L^2(\p D)$,
we see that trace $v|_{\p D}$ is well defined and is the
     constant function having
value $c_0$. Since restrictions of $h$ to $D$ and $\Omega\setminus
\overline D$ are in $H^1(D)$
and $H^1(\Omega\setminus \overline D)$, correspondingly,
and the trace from both sides of $\p D$ coincide,
we see that $h\in H^1(\Omega)$.

Next we show that $\sigma \nabla v\in L^2(\Omega\setminus \overline D)$.
Let $e_1,e_2,\dots,e_n$ be the standard Euclidean
    coordinate vectors in $\R^n$.
Then for $x\in \Omega\setminus \overline D$
\ba
\tilde g^{jk}\p_kv(x)=(e_j,\nabla_{\tilde g} v(x))_{\tilde g}=
(d(F^{-1})e_j,\nabla_g u(F^{-1}(x)))_{g}
\ea
and we see from (\ref{Q 4})
that this inner product is uniformly bounded.
By (\ref{cond and metr}) and (\ref{Q 4}),
$|\det(\tilde g)(x)|\leq c_3\,\hbox{dist}_{\R^n}\,(F^{-1}(x),y)^2$
and thus in $\Omega\setminus \overline D$ the functions
$V_j(x)=|\hbox{det}\, \tilde g(x)|^{1/2}
\delta_{jk} \tilde g^{ki}\p_i v$ satisfy
\ba
|V_j(x)|\leq c_4\,\hbox{dist}_{\R^n}\,(F^{-1}(x),y).
\ea
This implies that  $\vec V=(V_1(x),\dots,V_n(x))$
is in $H(\Omega\setminus \overline D;\nabla\cdotp)\cap
C(\overline \Omega\setminus D^{int})$. Moreover, by \cite{Gr}
the normal trace $\vec W\mapsto n\cdotp \vec W|_{\partial D}$
is a continuous map $H(\Omega\setminus \overline D;\nabla\cdotp)\cap
C(\overline \Omega\setminus D^{int})\to H^{-1/2}(\partial D)$,
we see that $n\cdotp \vec V|_{\partial D}=0$. Since the normal traces
of $\vec V$ coincide from both sides of $\partial D$
we see that $\vec V\in H(\Omega;\nabla\cdotp)$ and
\ba
\sum_{j,k=1}^n \p_k(|\hbox{det}\, \tilde g(x)|^{1/2}
\tilde g^{ki}\p_i h)=\nabla\cdotp \vec V=
0\quad \hbox{in }\Omega
\ea
in the sense of distributions.
This means that $h$ satisfies the  conductivity
equation in $\Omega$ in the sense of distributions.
This proves Theorem \ref{main}.
\medskip

A similar approach can be used to construct a different type of counterexample.
In contrast to the previous one, here the conductivity is not bounded
above near the
singular surface. Consider again the sets $\Omega=B(0,2)$
and $D=B(0,1)$ in $\R^3$.
In spherical coordinates
we define the metric $\hat g$ and the conductivity $\hat \sigma$
in the domain $\{1<r<2\}$
by the matrices
\ba
& & \hat g=
\left(\begin{array}{ccc}
(r-1)^{-2} & 0 & 0\\
0 & \rho^2 & 0 \\
0 & 0 &  \rho^2 \sin^2 \theta\\
\end{array}
\right),
\\
& &\hat \sigma=
\left(\begin{array}{ccc}
(r-1)\rho^{2}\sin \theta & 0 & 0\\
0 & (r-1)^{-1}\sin \theta & 0 \\
0 & 0 &  (r-1)^{-1}\sin^{-1} \theta\\
\end{array}
\right)
\ea
where $\rho>0$ is a constant.
Thus, $(\Omega\setminus \overline D,\hat g)$ is isometric to the product
$\R_+\times
S^2_\rho$ with the standard metric where $S^2_\rho$ is a 2-sphere with
radius $\rho$. Now  extend the conductivity $\hat \sigma$ to $\Omega$
so that it is positive definite in $D$.
It can be shown that, in the domain $\Omega\setminus \overline D$, the
equation
\ba
& &\nabla\cdotp \hat \sigma\nabla v(x)=0\hbox{ in }\Omega\setminus\overline D\\
& &v|_{\p \Omega}=f_0, \nonumber \\
& &v\in L^\infty(\Omega\setminus \overline D) \nonumber
\ea
has a unique solution. Also, when $\rho$ is small enough, we can
extend the definition of $v(x)$ to the whole domain $\Omega$
by defining it to have the constant value $c_0=\lim_{x\to \partial D} 
v(x)$ in $D$.
The  function $v(x)$ thus obtained
is a solution of the boundary value problem
\ba
& &\nabla\cdotp \hat \sigma\nabla v(x)=0\hbox{ in }\Omega,\\
& &v|_{\p \Omega}=f_0, \nonumber \\
& &v\in L^\infty(\Omega) \nonumber
\ea
in the sense of distributions.

In this case, the boundary measurements do not give information
about the metric inside $D$. To see this,  consider
$(\Omega\setminus \overline D,\hat g)$ as $\R_+\times S^2_\rho$
with coordinates $(t,\phi,\theta)$, where $t=t(r)=-\log (r-1)$, i.e.,
$t\in \R_+$ and
$\phi,\theta$ are coordinates on $S^2$.
We see that
for sets $B(0,1+R^{-1})\subset \Omega$, $R>1$
in original coordinates we can use
the superharmonic potentials
\ba
h(r,\phi,\theta)=\min(1,\frac {\log(r-1)}{\log(R-1)})
=\min(1,\frac {t(r)}{t(R)})
\ea
to see that the Brownian motion sent from $x\in \Omega\setminus
\overline D$ does not enter  $D$ a.s..
By writing the solution $u$
in terms of spherical harmonics, we have
\ba
u(t,\phi,\theta)=\sum_{n=1}^\infty \sum_{m=-n}^n a_{n,m}
e^{-\lambda_nt}Y^n_m(\phi,\theta)
\ea
where $\lambda_n=\rho^{-1}\sqrt{n(n+1)}$. Using
the original coordinates $(r,\phi,\theta)$ of $\Omega\setminus \overline D$
we see
\ba
u(r,\phi,\theta)=\sum_{n=1}^\infty \sum_{m=-n}^n a_{n,m}
(r-1)^{\lambda_n}Y^n_m(\phi,\theta).
\ea
Thus, if $\rho$ is small enough,  the solution
goes to constant and its derivative to zero
when $r\to 1^+$ so fast that  the solutions
of the conductivity equation  satisfy the equation  in
sense of distributions.

\bibliographystyle{amsalpha}

\bibliographystyle{amsalpha}

\end{document}